\documentclass[11pt]{amsart}
\usepackage{amsmath,paralist,amssymb}
\usepackage{amsthm}
\usepackage{hyperref}
\usepackage{graphicx,subfigure}
\usepackage{tikz,color}
\usepackage[margin=1in]{geometry}
\usetikzlibrary{arrows,shapes,snakes,backgrounds} 
\usepackage{verbatim}

\usepackage[all]{xy}

\newcommand{\old}[1]{}
\newtheorem{theorem}{Theorem}
\newtheorem*{corollary*}{Corollary}
\newtheorem*{theorem*}{Theorem}

\newtheorem{proposition}[theorem]{Proposition}
\newtheorem{lemma}[theorem]{Lemma}
\newtheorem{corollary}[theorem]{Corollary}
\newtheorem{conjecture}[theorem]{Conjecture}
 
\theoremstyle{remark}

\theoremstyle{definition}
\newtheorem{definition}{Definition}

\newcommand{\be}{\begin{equation}}

\newcommand{\ee}{\end{equation}}

\newcommand{\bd}{\begin{definition}}
\newcommand{\ed}{\end{definition}}
\newcommand{\bt}{\begin{theorem}}
\newcommand{\et}{\end{theorem}}
\newcommand{\bl}{\begin{lemma}}
\newcommand{\el}{\end{lemma}}
\newcommand{\bp}{\begin{proposition}}
\newcommand{\ep}{\end{proposition}}
\newcommand{\bc}{\begin{corollary}}
\newcommand{\ec}{\end{corollary}}

\newtheorem*{lemma*}{Lemma}

\newtheorem*{theorem1*}{Theorem \ref{sp}}
\newtheorem*{theorem2*}{Theorem \ref{mm}}

\def\i{{\it in}}
\def\out{{\it out}}

\def\In{{\it In}}
\def\Out{{\it Out}}

\def\fin{{\rm fin}}
\def\v{{\bf a}}

 \def\f_H{{\bf w}}
 \def\f{{\bf f}}
 \def\a{{\bf a}}
\def\R{\mathbb{R}}

\def\Z{\mathbb{Z}}

 \def\F{\mathcal{F}}
\def\O{\sigma}

\def\alt{{\rm dep}}

 \def\l{{\it i}}
 \def\dim{{\rm dim}}
 \def\x{{\bf x}}
 \def\a{\alpha}
 \def\b{\beta}
\def\R{\mathbb{R}}

\def\Z{\mathbb{Z}}
 \def\F{\mathcal{F}}
 
\def\C{\mathcal{C}}

\def\v{{\rm v}}
\def\ini{{\rm in }}
\def\fin{{\rm fin}}

 \def\f_H{{\bf w}}
 \def\f{{\bf f}}
 \def\a{\alpha}
\def\R{\mathbb{R}}

\def\Z{\mathbb{Z}}

 \def\F{\mathcal{F}}

\def\yb{\widetilde{ACYB_n}(\beta)}
\def\t{\mathcal{S}(\b)}

\begin{document}

\title[$h$-polynomials via  reduced forms]{$h$-polynomials via reduced forms}
\author{Karola M\'esz\'aros}
\email{karola@math.cornell.edu}
\address{
Department of Mathematics, Cornell University, Ithaca, NY 14853
}
\thanks{The author was partially supported by a National Science Foundation Postdoctoral Research Fellowship  (DMS 1103933)}
\date{April 7, 2015}

\begin{abstract} The flow polytope  $\F_{\widetilde{G}}$ is the set of nonnegative unit flows on the graph $\widetilde{G}$. The subdivision algebra of flow polytopes prescribes a way to dissect a flow polytope $\F_{\widetilde{G}}$  into simplices. Such a dissection is encoded by the terms of the so called reduced form of the monomial $\prod_{(i,j)\in E(G)}x_{ij}$. We prove that we can use the subdivision algebra of flow polytopes to construct not only dissections, but also regular flag triangulations of flow polytopes. We prove that reduced forms in the subdivision algebra are generalizations of $h$-polynomials of the triangulations of flow polytopes. We deduce several corollaries of the above results, most notably proving certain cases of a conjecture of Kirillov about the nonnegativity of reduced forms in the noncommutative quasi-classical Yang-Baxter algebra. 


\end{abstract}

 \maketitle
 \tableofcontents 
 
 \section{Introduction}
 \label{sec:intro}
 
 Nonnegativity properties abound in mathematics, and whenever one arises, the most satisfying explanation of integer nonnegativity is to demonstrate what a certain nonnegative quantity  counts. The present paper is written in this spirit and explains nonnegativity properties of polynomials using geometric interpretations for their coefficients. 
 
 We study dissections of flow polytopes via the subdivision algebra, and show that we can construct  regular flag triangulations of flow polytopes using the subdivision algebra. This in turn empowers us to prove several interesting corollaries, one in particular partially proving a nonnegativity conjecture of Kirillov \cite[Conjecture 2]{k2} of certain polynomials  called reduced forms in the quasi-classical Yang-Baxter algebra. 
 
With our geometric methods  we can study reduced forms  in the subdivision algebra and the quasi-classical Yang-Baxter algebra. The latter algebra was introduced by A.N. Kirillov \cite{Kir, k2, k2014} with Schubert calculus in mind and it  is closely related to the Fomin-Kirillov algebra \cite{fk}.  The subdivision algebra has been considered by the present author under this name, as its relations encode ways to subdivide root and flow polytopes \cite{m1, m2, m-prod}. It has also been considered by Kirillov \cite{Kir, k2, k2014}  since it is the abelianization of the  quasi-classical Yang-Baxter algebra (and this is how he refers to it). The polynomials of interest in this paper arise as reduced forms in the above algebras; the reduced form 
 of a monomial in an algebra is obtained via substitution rules dictated by the relations of the algebra.

 The essence of the subdivision algebra is that the reduced form of a monomial in it can naturally be seen as a dissection of a flow polytope corresponding to the monomial into simplices. We show that among these dissections that {the reduced form can encode} are also   {unimodular, regular and flag triangulations of flow polytopes}, see Theorem \ref{tri}. Using the connection of reduced forms to dissections 
 we show that {reduced forms are multivariate generalizations of $h$-polynomials}, see Theorem \ref{thm-h}.  
 Recall that the $h$-polynomial of a simplicial complex is a way of encoding the number of faces of each dimension.  We prove that   if we set certain variables of  reduced forms to $1$  in the subdivision algebra  we obtain the  (shifted) $h$-polynomials of regular triangulations of flow polytopes.  This result opens a new avenue for understanding reduced forms of monomials in the subdivision and related algebras.     
 We prove nonnegativity results in the subdivision algebra as a consequence of the specialized reduced form equaling the shifted $h$-polynomial, see Theorem \ref{pos}.  As a corollary  we establish a special case of Conjecture 2 of Kirillov appearing in \cite{k2} about the nonnegativity of reduced forms in  the quasi-classical Yang-Baxter algebra, see Theorem \ref{special}. We also express specialized reduced forms in terms of Ehrhart series of flow polytopes, which in turn can be seen in terms of Kostant partition functions, see Theorem \ref{thm-shell2} and Lemma \ref{thm-ehr}.

 Our methods in this paper are largely geometric. In  \cite{h-poly2} we study reduced forms from the point of view of the structure of reduction trees, leaving the geometry behind.

  The paper is organized as follows. In Section \ref{sec:def} we define flow polytopes. Next we explain how to subdivide flow polytopes and how we can encode the subdivisions with a reduction tree. Then we define the  subdivision algebra and show that the reduced form can be read off from the leaves of the  reduction tree.

   In Section \ref{sec:tri} we show that we can use the subdivision algebra and arrive not only at a dissection of a flow polytope, but also to a triangulation of the flow polytope, in the sense of a simplicial complex. To do this we use a particular reduction order $\O$. The triangulation we obtain is regular and flag. In Section \ref{alt}    we  prove that the reduced form of a monomial in the subdivision algebra specialized at certain variables is equal to the shifted $h$-polynomial of  the aforementioned triangulation of the flow polytope of the associated to the monomial.  
      In Section \ref{sec:leaves} we   describe the full set of leaves of the reduction tree in order $\O$, or equivalently, the monomials in the reduced form of  a monomial. In
 Section \ref{sec:nonneg} we use the regularity of the triangulation we constructed to prove and interpret the nonnegativity of the coefficients of the reduced form of a monomial in the subdivision algebra as well as a special case of Conjecture 2 of Kirillov appearing in \cite{k2}.  In Section \ref{sec:ehr}   we relate reduced forms to Ehrhart series of flow polytopes, and thus  obtain a generalization of  \cite[Theorem 3.10]{k2}. We also relate reduced forms to Kostant partition functions.

 \section{Definitions and more of the story}
\label{sec:def}

\subsection{Flow polytopes and their subdivisions.} Given a loopless graph $G$  on the vertex set $[n]$,  let $\ini(e)$ denote the smallest (initial) vertex of edge $e$ and $\fin(e)$ the biggest (final) vertex of edge $e$.  Let $E(G)=\{\{e_1, \ldots, e_l\}\}$   be the multiset of edges of $G$. We correspond  variables $x_{e_i}$, $i \in [l]$, to the edges of $G$, of which we think as flows.  
The {\bf flow polytope} $\F_G$ is naturally embedded into $\R^l$, where  $x_{e_i}$, $i \in [l]$, are thought of as the coordinates. $\F_G$  is defined by 

 $$x_{e_i}\geq 0, \mbox{ }i \in [l],$$

  $$1=\sum_{e \in E(G), \ini(e)=1}x_e= \sum_{e \in E(G),  \fin(e)=n}x_e,$$
  
\noindent  and for $2\leq i\leq n$
  
  $$\sum_{e \in E(G), \fin(e)=i}x_e= \sum_{e \in E(G), \ini(e)=i}x_e.$$

 The flow polytope $\F_{K_{n+1}}$, where $K_{n+1}$ is the complete graph on $n+1$ vertices,  can be thought of as the Chan-Robbins-Yuen polytope \cite{CRY}, and  has received a lot of attention, since  its volume is equal to $\prod_{k=0}^{n-2}Cat(k),$ where  $Cat(k)=\frac{1}{k+1}{2k \choose k}$ is the $k$th Catalan number. There  is no combinatorial proof of the aforementioned result;  Zeilberger \cite{Z} provided an analytical proof. For more of the story see \cite{mm}.

Flow polytopes lend themselves to subdivisions via  {\it reductions}, as explained below. A similar property of root polytopes was studied in \cite{m2, m1}.

\bd \label{3} Given a graph $G$ on the vertex set $[n]$  containing edges $(i, j)$ and $(j, k)$, $i<j<k$, performing the {\bf reduction} on these edges of $G$ yields three graphs on the vertex set $[n]$:
  \begin{eqnarray} \label{graphs1}
E(G_1)&=&E(G)\backslash \{(j, k)\} \cup \{(i, k)\}, \nonumber \\
E(G_2)&=&E(G)\backslash \{(i, j)\} \cup \{(i, k)\},  \nonumber \\
E(G_3)&=&E(G)\backslash \{(i, j), (j, k)\} \cup \{(i, k)\}.
\end{eqnarray}

 When performing a reduction on the edges $(i, j), (j, k)$  we say that the edge $(i,j)$ is {\bf dropped} if we go towards $G_2$ or $G_3$ as in \eqref{graphs1} and  $(i,j)$ is {\bf kept} if we go towards $G_1$. Similarly, edge $(j,k)$ is dropped if we go towards $G_1$ or $G_3$ as in \eqref{graphs1} and  $(j,k)$ is  kept if we go towards $G_2$.

\ed

\begin{figure}
\begin{center}
\includegraphics[scale=.7]{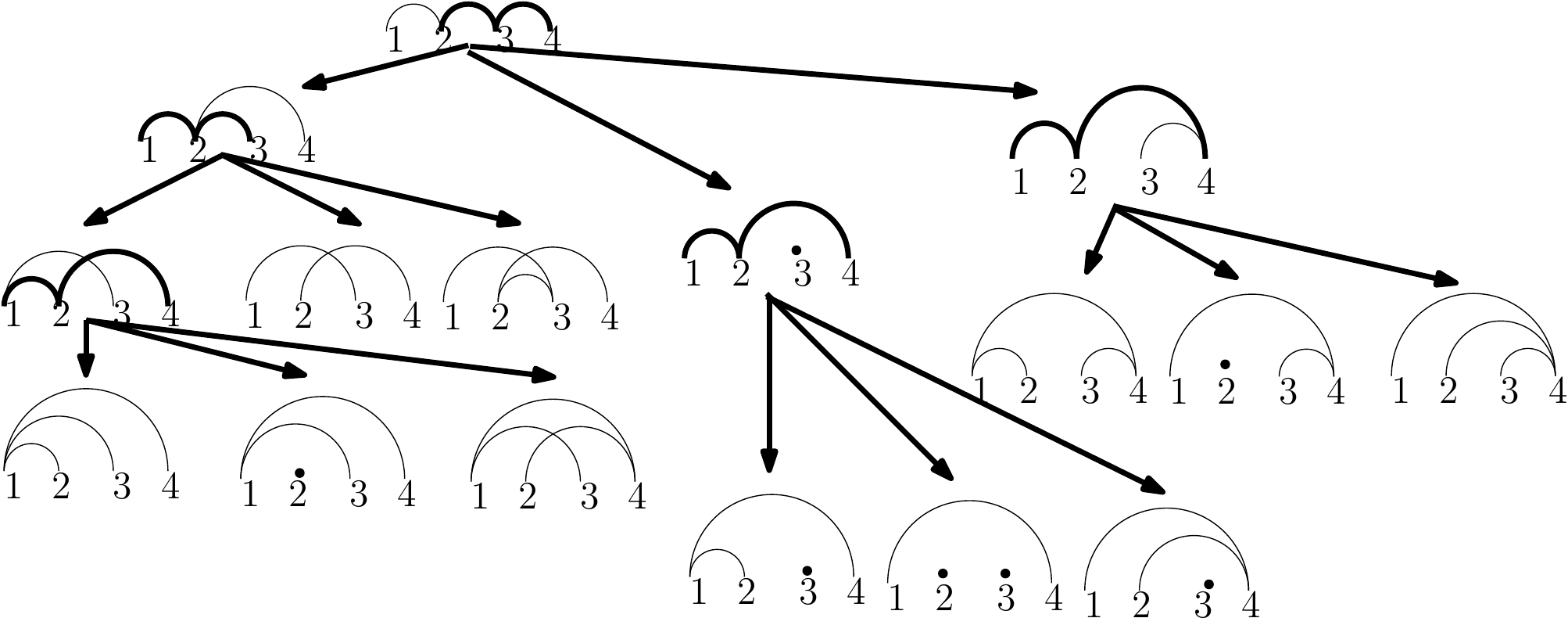}
 \caption{A reduction tree of $G=([4], \{(1,2),(2,3), (3,4)\})$. The edges on which the reductions are performed are in bold. }
 \label{fig:redtree}
 \end{center}
\end{figure}

\bd A \textbf{reduction tree} $R_G$ of a graph $G$ is a tree with nodes labeled by graphs and such that all non-leaf nodes of $R_G$ have three children. The root is labeled by $G$. If there are two edges  $(i, j), (j, k) \in E(G)$,  $i<j<k$, on which we choose to do a reduction, then the  children of the root are labeled by  $G_1, G_2$ and $G_3$ as in (\ref{graphs1}). Next, continue this way by constructing  reduction trees for $G_1, G_2$ and $G_3$. If some graph has no edges $(i, j), (j, k)$,  $i<j<k$, then it is its own reduction tree. Note that the reduction tree $R_G$ is not unique; it depends on our  choice of edges to reduce.
However, the number of leaves (referring to the graph labeling a leaf) of all reduction trees of $G$ with a given number of edges is the same as Lemma \ref{leaves} states below. We choose a particular embedding of the reduction tree in the plane for convenience: we root it at $G$ with the tree growing downwards, and such that the left child is $G_1$, the middle child is $G_3$ and the right child is $G_2$; see Figure \ref{fig:redtree}. The leaves which have the same number of edges at the root are called \textbf{full dimensional}.
\ed

\bd \label{edges}
Let the edges of $G$ be $e_1, \ldots, e_k$, where we distinguish multiple edges. If a reduction involving edges $a=(i,j)$ and $b=(j,k)$ of $G$ is performed, then the new edge $(i,k)$ appearing in all three graphs as in \eqref{graphs1} is formally thought of as $a+b(=b+a)$. The other edges stay unchanged. To get to nodes $G_1$ and $G_2$  of $R_G$ we iterate this process, thereby expressing the edges of any node as a sum of edges of the graph being the root of the reduction tree. Two edges $c$ and $d$ in the graphs  $G_1$ and $G_2$, respectively, are the same, if they are the sum of exactly the same edges of $G$. The intersection of two graphs $G_1$ and $G_2$ in a reduction tree $R_G$ is $G_1\cap G_2=(V(G), E(G_1)\cap E(G_2))$, where if $e \in E(G_1)\cap E(G_2)$ then as explained above $e$ is the sum of the same edges of $G$ in both $G_1$ and $G_2$. 
\ed

\bl \label{leaves} \cite{m2} Given two distinct reduction trees of $G$, let $r^1_k$ and $r^2_k$ be the number of leaves with $k$ edges is them, respectively. Then, $r^1_k=r^2_k$. 
  \el
 
 \old{
To prove Lemma \ref{leaves}, recall that given a  polytope $\mathcal{P}\subset \mathbb{R}^{N}$, the {$t^{th}$ dilate} of $\mathcal{P}$ is $\displaystyle t \mathcal{P}=\{(tx_1, \ldots, tx_{N}) \mid  (x_1, \ldots, x_{N}) \in \mathcal{P}\}.$ The number of lattice points of $t\mathcal{P}$, where $t$ is a nonnegative integer and $\mathcal{P}$ is a convex polytope, is given by the {\bf Ehrhart function} $i({\mathcal{P}}, t)$. If $\mathcal{P}$ has integral vertices then $i({\mathcal{P}}, t)$ is a  polynomial.

\medskip

 \noindent {\it Proof of  Lemma \ref{leaves}.} Let $f_k$ denote the number of leaves of a reduction tree of $G$ which are graphs with $k$ edges. Then, 
 $$i(\F_{{\widetilde{G}}}^\circ, t)=\sum_{k} f_k {t-1 \choose n+4-k},$$ since by \cite{BV2} the dimension of a flow polytope is $\#V-\#E+1$, and the Ehrhart polynomial of the $t$th dilate of the interior of a $d$-dimensional simplex corresponding to a graph labeling the  leaf of a reduction tree is ${t-1 \choose d}$ \cite{BR}. But since the binomial coefficients $ {t-1 \choose n+4-k},$ $n+4-k=0, 1, \ldots$  form a basis of $\Z[t]$, it follows that $f_k$ is unique and independent of the particular reduction tree. 
  \qed

}
 
\bd  The {\bf augmented graph} ${\widetilde{G}}$ of $G=([n], E)$ is ${\widetilde{G}}=([n]\cup \{s, t\}, \widetilde{E})$, where $s$ (source) is the smallest, $t$ (target/sink) is the biggest vertex of  $[n]\cup \{s, t\}$, and $\widetilde{E}=E\cup \{(s, i), (i, t) | i \in [n]\}$.  Denote by $\mathcal{P}(\widetilde{G})$ the set of all maximal paths in $\widetilde{G}$, referred to as \textbf{routes}. It is well known that the unit flows sent along the routes in $\mathcal{P}(\widetilde{G})$ are the vertices of $\mathcal{F}(\widetilde{G})$. 
 \ed

\begin{definition} \label{leaf} Consider a node $G_1$ of the reduction tree $R_G$, where each edge of $G_1$ is considered as a sum of the edges of $G$. The image of the map $m:E(G_1)\rightarrow \mathcal{P}(\widetilde{G})$ which takes an edge $(v_1, v_2)=e=e_{i_1}+\cdots+e_{i_l}$, $e \in G_1$, $e_{i_j} \in E(G)$, $j \in [l]$, to the route $(s, v_1), e_{i_1},\ldots,e_{i_l}, (v_2, t)$ gives the vertices of $\F_{\widetilde{G_1}}$ (by taking the unit flows  on these routes).  In case $G_1$ is not a node of the reduction tree $R_G$, but it is  an intersection of  nodes of $R_G$, so that each edge of $G_1$ can still be considered as a sum of the edges of $G$, we still define $\F_{\widetilde{G_1}}$ as above. This definition of  $\F_{\widetilde{G_1}}$ is of course with respect to $G$, and this is understood from the context.
\end{definition}

Using the above definitions the proof of  the following lemma is an easy exercise.

\bl{ \cite[Proposition 1]{m-prod},\cite[Proposition 4.1]{mm}, \cite{p, S}}
 \label{red} Given a graph $G$ on the vertex set $[n]$ and   $(i, j), (j, k) \in E(G)$,   for some $i<j<k$,  and $G_1, G_2, G_3$ as in \eqref{graphs1} and $ \F_{{\widetilde{G}_i}}$, $i \in [3]$, as in Definition \ref{leaf} we have  $$\F_{{\widetilde{G}}}=\F_{{\widetilde{G}_1}} \bigcup \F_{{\widetilde{G}_2}},  \F_{{\widetilde{G}_1}} \bigcap \F_{{\widetilde{G}_2}}=\F_{{\widetilde{G}_3}} \text{ and } \F_{{\widetilde{G}_1}}^\circ \bigcap \F_{{\widetilde{G}_2}}^\circ=\emptyset,$$

 \noindent where $\F_{{\widetilde{G}}}$, $\F_{{\widetilde{G}_1}}$, $\F_{{\widetilde{G}_2}}$ are of the same dimension $d-1$, $\F_{{\widetilde{G}_3}}$ is $d-2$ dimensional,  and  $\mathcal{P}^\circ$ denotes the interior of $\mathcal{P}$.
\el

 
\subsection{Encoding subdivisions by relations.} Note that the reduction of graphs given in \eqref{graphs1} can be encoded as the following relation:

\be x_{ij}x_{jk}=x_{ik}x_{ij}+x_{jk}x_{ik}+\b x_{ik},  \mbox{ for }1\leq i<j<k\leq n. \ee

Namely, interpreting the double indices of the variables $x_{ij}$ as edges, the monomial $x_{ij}x_{jk}$ picks out two edges $(i,j), (j, k)$, $i<j<k$, and replaces it with three monomials, corresponding to operation on graphs  \eqref{graphs1}. The variable $\beta$ is simply a placeholder, indicating that the number of edges in the third graph is one less than in the other graphs

These relations give rise to what we call the {\it subdivision algebra}. 

\bd The associative \textbf{subdivision algebra}, denoted by $\t$, is an associative algebra, over the ring of polynomials $\Z[\b]$, generated by the set of elements $\{x_{ij} : 1 \leq i<j \leq n\}$, subject to the relations:

(a) $x_{ij}x_{kl}=x_{kl}x_{ij},$ if  $i<j$, $k<l$,

(b) $x_{ij}x_{jk}=x_{ik}x_{ij}+x_{jk}x_{ik}+\b x_{ik},$ if $1\leq i<j<k\leq n$.

\ed

The algebra $\t$ has been studied  in the context of root polytopes \cite{m2}.

\bd \label{rf} Given a monomial $M$  in $\t$, its reduced form is defined as follows.  Starting with $p_0=M$, produce a sequence of polynomials $p_0, p_1, \ldots, p_m$ in the following fashion.  To obtain $p_{r+1}$ from $p_r$,  choose a monomial of   $p_r$ which  is divisible by $x_{ij}x_{jk}$, for some $i,j,k$, and replace the factor  $x_{ij}x_{jk}$ in this monomial   with   $x_{ik}x_{ij}+x_{jk}x_{ik}+\b x_{ik}$. Note that   $p_{r+1}$   has two more monomials than $p_r$. Continue  this process until   a  polynomial  $p_m$  is obtained, in which  no monomial is divisible by  $x_{ij}x_{jk}$, for any $i,j,k$.  Such a polynomial $p_m$  is a {\bf reduced form} of $M$. Note that we allow the use of the commutation relations  in this process.
\ed

Given a monomial $M$ we can encode it by a graph $G_M$, simply by letting the edges of $G_M$ be the given by the indices of the variables in $M$. Denote a reduced form of $M$  in $\t$ by $Q_{G_M}^{\t}(\x; \b)$. Note that the reduced form of $M$ is not unique, so if we wanted to specify the reduced form exactly, we would also need to specify a reduction tree for $G$. When writing $Q_{G_M}^{\t}(\x; \b)$ we pick an arbitrary reduced form at hand.  If in the reduced forms we set $\x=(1, \ldots, 1)$, then in the notation we omit $\x$:  
$Q_{G_M}^{\t}(\b)$.  We will show later that $Q_{G_M}^{\t}(\b)$ is independent of the choice of reduction tree, thus $Q_{G_M}^{\t}(\b)$ is well-defined as is.

It is easy to see that by definition, a reduced form of a monomial in the subdivision algebra can be read off from the leaves of the reduction tree of the corresponding graph. With this and Lemma \ref{red} in mind, it is no surprise that we can  prove results about reduced forms of monomials in the algebras using flow polytopes.

Recall  that a reduced form of a monomial in  $\t$ is not necessarily unique, which could be a desirable property. Amazingly, there is  a noncommutative algebra, denoted $\yb$, which is much like $\t$, yet in which the reduced forms are unique. The latter statement was proved  in \cite{m2}.  It was A.N. Kirillov \cite{k2, k2014} who introduced $\yb$ and shed the first light on its rich combinatorial structure. The paper \cite{h-poly2}  addresses more of the story of subdivision algebras as well as the story of  $\yb$.

 \section{Triangulating flow polytopes}
 \label{sec:tri}

In this section we prove  that we can use the subdivision algebra to obtain triangulations of every flow polytope $\F_{\widetilde{G}}$. A priori this is far from clear -- we are only guaranteed dissections of flow polytopes, but not the simplicial complex structure.  Moreover, the triangulation we obtain is flag and regular, and thus shellable. 
We construct our special triangulation by picking a specific reduction order $\O$ on our graph $G$ and utilizing the properties of $\O$ to prove properties about the leaves of the reduction tree $R_G^{\O}$ obtained using $\O$.  We then note that there is a whole class of orders for which our arguments work, but for clarity we lay out the argument for $\O$ first. 


\subsection{A family of orders $\F(\O)$.}
\label{sec:orders}

\begin{figure}
\begin{center}
\includegraphics[scale=.9]{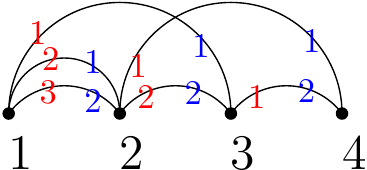}
 \caption{The ordering of the incoming and outgoing edges at each vertex from topmost to bottom.}
 \label{fig:sigma}
 \end{center}
\end{figure}

  Given an arbitrary graph $G$ on the vertex set $[n]$, put a total order on the set of incoming and as well as a total order set of outgoing edges at each vertex $1<v<n$ (thus we consider each edge twice here, once as an incoming edge and once as an outgoing edge).
 Do the reductions in $G$ proceeding from the smallest vertex towards the greatest in order. Look for the smallest vertex $v$ which is nonalternating, that is that has both an edge $(a, v)$ and an edge $(v, b)$ incident to it, with $a<v<b$. (Two edges $(a, v)$ and   $(v, b)$ are \textbf{non-alternating} if $a<v<b$.) Look at the incoming  and outgoing edges at $v$,  $(a, v)$  and  $(v, b)$, which are smallest in the ordering of the incoming and outgoing edges, respectively. Do the reduction on $(a, v)$  and  $(v, b)$. In the three obtained graphs the relative ordering of the edges stays the same, with the new edge $(a, b)$ either taking the place of   $(a, v)$  or  $(v, b)$ if these were dropped, or directly preceeding them when they are kept.  Continue in this fashion on each leaf of the partial reduction tree ultimately arriving to the reduction tree $R_G$ with all leaves alternating graphs, that is all of their vertices are alternating.  Let $\O$ be the order where the initial ordering of the incoming and  outgoing edges at each vertex is such that the topmost is the smallest, then the next topmost, etc. See Figure \ref{fig:sigma} for an example of this ordering. All results of this paper generalize for any order in $\F(\O)$, but for simplicity we state and prove them for $\O$ only.

 \old{
\begin{figure}
\begin{center}
\includegraphics[scale=.7]{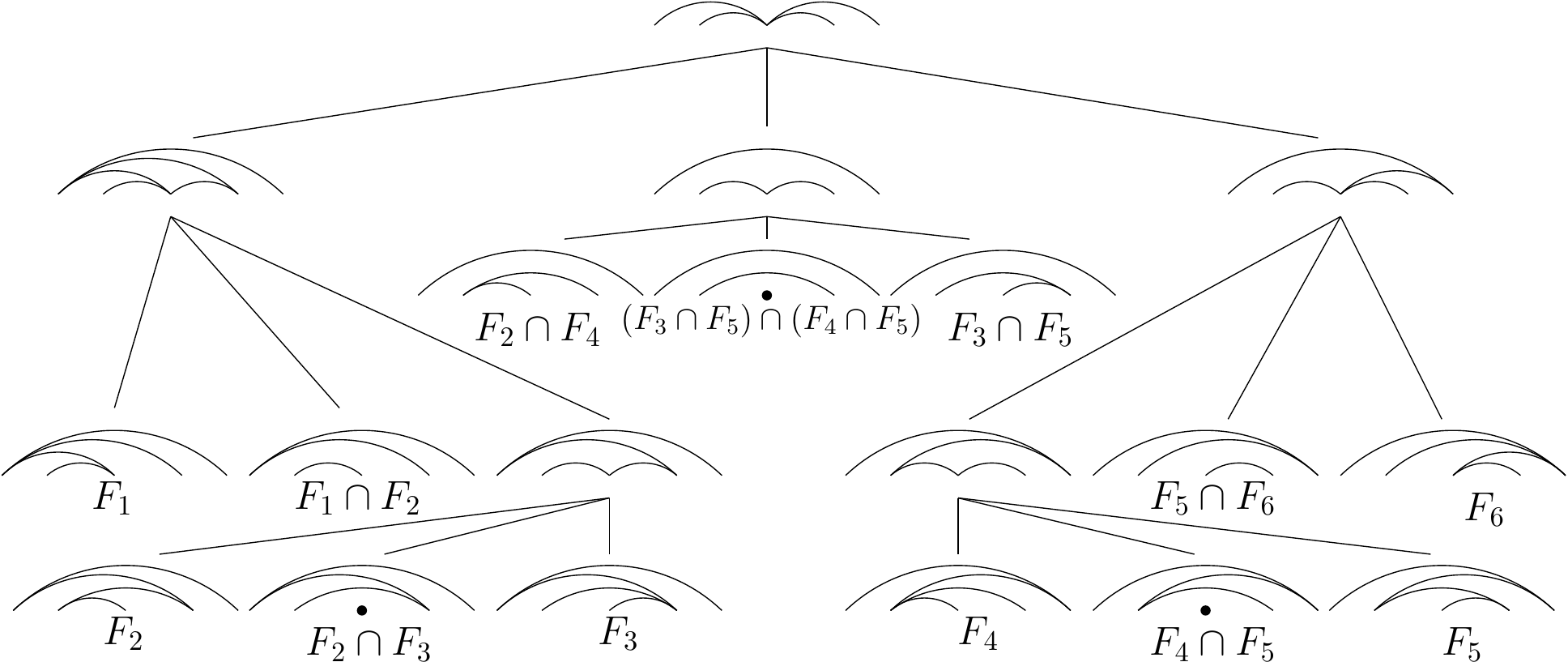}
 \caption{The reduction tree of $G=([5], \{(1,3),(2,3), (3,4), (3,5)\})$ with reductions executed in order $\O$.  The labels $F_i$, $i \in [6]$,  are explained in Theorems \ref{thm-shell2} and \ref{thm-shell}.}
 \label{fig:o}
 \end{center}
\end{figure}
}

The main result of this section is:

\begin{theorem} \label{tri} The simplices corresponding to the full dimensional leaves of $R_G^{\O}$  yield the top dimensional simplices in a regular and flag triangulation of $\F_{\widetilde{G}}$. Moreover, lower dimensional simplices of this triangulation which are not contained in the boundary of $\F_{\widetilde{G}}$ are obtained from the (not  full dimensional) leaves of $R_G^{\O}$.
\end{theorem}

\proof The proof follows from Theorem \ref{k}, Lemma \ref{flag} and Proposition \ref{coh}.
\qed

\medskip

Next we explain Theorem \ref{k} and prove Proposition \ref{coh}.  

\old{
\bt \label{coord} Let $G_1$ and $G_2$ be two leaves of $R_G^{\O}$. Then 

\be \label{int}  \F_{\widetilde{G_1}}\cap \F_{\widetilde{G_2}}=\F_{\widetilde{G_1\cap G_2}}.\ee   \et

Note that $\F_{\widetilde{G_1\cap G_2}} \subset \F_{\widetilde{G_1}}\cap \F_{\widetilde{G_2}}$ follows for any reduction order by definition of these polytopes. However, the other direction is not clear, since an $R_G$-simplicial collection need not be a simplicial complex, and so the simplices do not have to intersect in their faces. We use results of \cite{}
}

\subsection{Coherent routes, cliques and triangulation of flow polytopes.} Given a graph $G$ on the vertex set $[n]$ with edges oriented from smaller to bigger vertices, the vertices of the flow polytope $\F_G$ correspond to integer flows  of size one on maximal directed paths from the source ($1$) to the sink ($n$). Following \cite{kosh} we call such maximal paths \textbf{routes}. The following definitions  follow \cite{kosh}. Fix a \textbf{framing} at each inner vertex $v$ (that is a vertex that is not a source or a sink) of $G$, which is the linear ordering $\prec_{\i(v)}$ on the set of incoming edges $\i(v)$ to $v$  and the linear ordering $\prec_{\out(v)}$ on the set of  outgoing edges $\out(v)$ from $v$.  We call a graph with a framing at each inner vertex framed. For a framed graph $G$ and an inner vertex $v$ we denote by $\In(v)$ and by $\Out(v)$ the set of maximal paths ending in $v$ and   the set of maximal paths starting at $v$, respectively. We define the order $\prec_{\In(v)}$ on the paths in $\In(v)$ as follows. If $P, Q \in \In(v)$  then let $w$ be the largest vertex after which $P$ and $Q$ coincide and before which they differ. Let $e_P$ be the edge of $P$ entering $w$ and $e_Q$ be the edge of $Q$ entering $w$. Then $P \prec_{\In(v)} Q$ if and only if $e_P \prec_{\i(w)} e_Q$. The linear order $\prec_{\Out(v)}$ on $\Out(v)$ is defined analogously. 

Given a route $P$ with an inner vertex $v$ denote by $Pv$ the maximal subpath of $P$ entering $v$ and by $vP$ the maximal subpath of $P$ leaving $v$. We say that the routes $P$ and $Q$ are \textbf{coherent at a vertex} $v$ which is an inner vertex of both $P$ and $Q$ if the paths $Pv, Qv$ are ordered the same way as $vP, vQ$; e.g., if 
$Pv \prec_{\In(v)}Qv$ and 
$vP \prec_{\Out(v)}vQ$. 
We say that routes $P$ and $Q$ are \textbf{coherent} if they are coherent at each common inner vertex. We call a set of mutually coherent routes a \textbf{clique}. The following theorem is a special case of \cite[Theorems 1 \& 2]{kosh}.

\begin{theorem}  \cite[Theorems 1 \& 2]{kosh} \label{k} Given a framed graph $G$, taking the convex hulls of the vertices corresponding to the routes in maximal cliques yield the top dimensional simplices in a regular triangulation of $\F_G$. Moreover, lower dimensional simplices of this triangulation are obtained as convex hulls of the vertices corresponding to the routes in (not  maximal) cliques.
\end{theorem}

\bl \label{flag} The triangulation described in Theorem \ref{k} is flag.
\el

\proof Consider a  non-face $N$ of the triangulation that is of cardinality greater than $2$. By Theorem \ref{k} the vertices of $N$ are routes that are not coherent; in particular there are two routes $P$ and $Q$ which yield vertices of $N$ and are not coherent. Since $P$ and $Q$ are not coherent, they constitute a non-face. Therefore, all minimal non-faces of the triangulation described in Theorem \ref{k} are of size $2$, and therefore the triangulation is flag.
\qed

\begin{figure}
\begin{center}
\includegraphics[scale=.9]{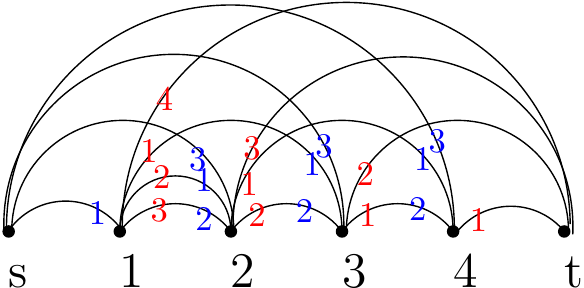}
 \caption{Framing $\widetilde{\O}$ on $\widetilde{G}$ with $G=([4], \{\{(1,2), (1,2), (1,3), (2,3), (2,4), (3,4)\}\})$.}
 \label{fig:frame}
 \end{center}
\end{figure}

Now we are ready to prove the proposition which together with Theorem \ref{k} and Lemma \ref{flag} implies Theorem \ref{tri}. We define the framing $\widetilde{\O}$ on $\widetilde{G}$ as the ordering $\O$ on the edges of $G$, that is, the incoming edges are ordered top to bottom and the outgoing edges are also ordered top to bottom, and the edges of the form $(s, i)$ and $(i, t)$, for $i \in [n]$, are always last in the orderings. See Figure \ref{fig:frame} for an example.

\begin{proposition} \label{coh} The set of vertices of the simplices corresponding to the leaves of $R_G^{\O}$ form a clique of mutually coherent paths in $\widetilde{G}$ with the framing $\widetilde{\O}$. 
\end{proposition}

\proof Suppose that to the contrary, there are two vertices of a simplex corresponding to a leaf of $R_G^{\O}$, which correspond to non-coherent routes $P$ and $Q$ in $\widetilde{G}$. Suppose that $P$ and $Q$ are not coherent at the common inner vertex $v$. Suppose that the smallest vertex after which $Pv$ and $Qv$ agree is $w_1$ and the largest vertex before which $vP$ and $vQ$ agree is $w_2$. Let the edges incoming to $w_1$ be $e_P^1$ and $e_Q^1$ for $P$ and $Q$, respectively, and let  the edges outgoing from $w_2$ be $e_P^2$ and $e_Q^2$ for $P$ and $Q$, respectively. Since $P$ and $Q$ are not coherent at $v$, this implies that either $e_P^1\prec_{\i(w_1)} e_Q^1$ and $e_Q^2\prec_{\out(w_2)} e_P^2$ or $e_Q^1\prec_{\i(w_1)} e_P^1$ and $e_P^2\prec_{\out(w_2)} e_Q^2$. We also have that the segments of $P$ and $Q$ between $w_1$ and $w_2$ coincide. Note that since the edges of the form $(s, i)$ and $(i, t)$, $i \in [n]$, are last in the linear orderings of the incoming and outgoing edges, it follows that at most one of the edges $e_P^1$ and $e_P^2$ and at most one of the edges $e_Q^1$ and $e_Q^2$ could be incident to $s$ or $t$. We consider several cases based on whether any of $e_P^1, e_P^2, e_Q^1, e_Q^2$ are incident to $s$ or $t$. Denote by $p$ the sum of edges between $w_1$ and $w_2$ on $P$. If none of $e_P^1, e_P^2, e_Q^1, e_Q^2$ are incident to $s$ or $t$, then after a certain number of reductions executed according to $\O$ we are about the perform the reduction $(*(e_{\bar{Z}}^1+p), e_Z^2)$, where $*(e_{\bar{Z}}^1+p)$ denotes the sum of edges left of $w_1$ that are edges in $\bar{Z}$ not incident to $s$ (including $e_{\bar{Z}}^1$ in particular) and  $p$,  $\{\bar{Z}, Z\}=\{P, Q\}$. Note, however, that after executing this reduction we have to drop either the edge $*(e_{\bar{Z}}^1+p)$ or the edge $e_Z^2$. However, if  the former were true, it would make it impossible for $*(e_{\bar{Z}}^1+p)+e_{\bar{Z}}^2$ to be a subsum in an edge of the leaf we are considering, which it has to be in order for $\bar{Z}$ to be a route giving a vertex of the simplex we are considering. The latter on the other hand would make it impossible for $*(e_{Z}^1+p)+e_{{Z}}^2$ to be a subsum in an edge of the leaf we are considering, where  $*(e_{{Z}}^1+p)$ denotes the sum of edges left of $w_1$ that are edges in ${Z}$ not incident to $s$. However then  ${Z}$ cannot be a route giving a vertex of the simplex we are considering. Thus we see that $Z$ and $\bar{Z}$, aka, $P$ and $Q$, cannot be incoherent in this way. It follows that we need to consider the possibilities where some of $e_P^1, e_P^2, e_Q^1, e_Q^2$ are incident to $s$ or $t$. One can construct  similar arguments to the above in all those cases. 
\qed

\section{Reduced forms are $h$-polynomials} 
\label{alt} 

In this section we show that for a graph $G$ the reduced form of the monomial $m[G]=\prod_{(i,j) \in E(G)}x_{ij}$ can be seen as the shifted $h$-polynomial of a unimodular triangulation of the flow polytope  $\F_{\widetilde{G}}$.

\bd Let $\C$ be a $d-1$ dimensional simplicial complex.  The \textbf{$f$-vector} of $\C$ is $$f(\C)=(f_{-1}, f_0, \ldots, f_{d-1}),$$ where $f_i=f_i(\C)$ be the number of $i$-dimensional simplices in $\C$. By convention, $f_{-1}=1$ unless $\C=\emptyset$, in which case $f_{-1}=0$. The \textbf{$h$-vector} of $\C$ is $h(\C)=(h_0, h_1, \ldots, h_d)$, defined by 

\begin{equation} \label{f-h} \sum_{i=0}^d f_{i-1}(x-1)^{d-i}=\sum_{i=0}^d h_{i}x^{d-i}.\end{equation}

Define the \textbf{$h$-polynomial} of a simplicial complex $\C$ to be 

\begin{equation} h(\C, \x)=\sum_{i=0}^d h_{i}x^{i}.
\end{equation}

\ed

\bl \label{l-h-poly} \cite{Stcom} For a simplicial complex $\C$ we have:
\begin{equation} \label{h-poly} h(\C, \x)=\sum_{F \in \C}x^{\#F}(1-x)^{d-\#F}.
\end{equation}
\el

The main result of this section is the following theorem. 

\begin{theorem} \label{thm-h} We have 
\be Q_G^{\t}(\b)=h(\C, \b+1),  \ee where $\C$ is any unimodular triangulation of ${\F_{\widetilde{G}}}$.
\end{theorem}

\begin{figure}
\begin{center}
\includegraphics[scale=.8]{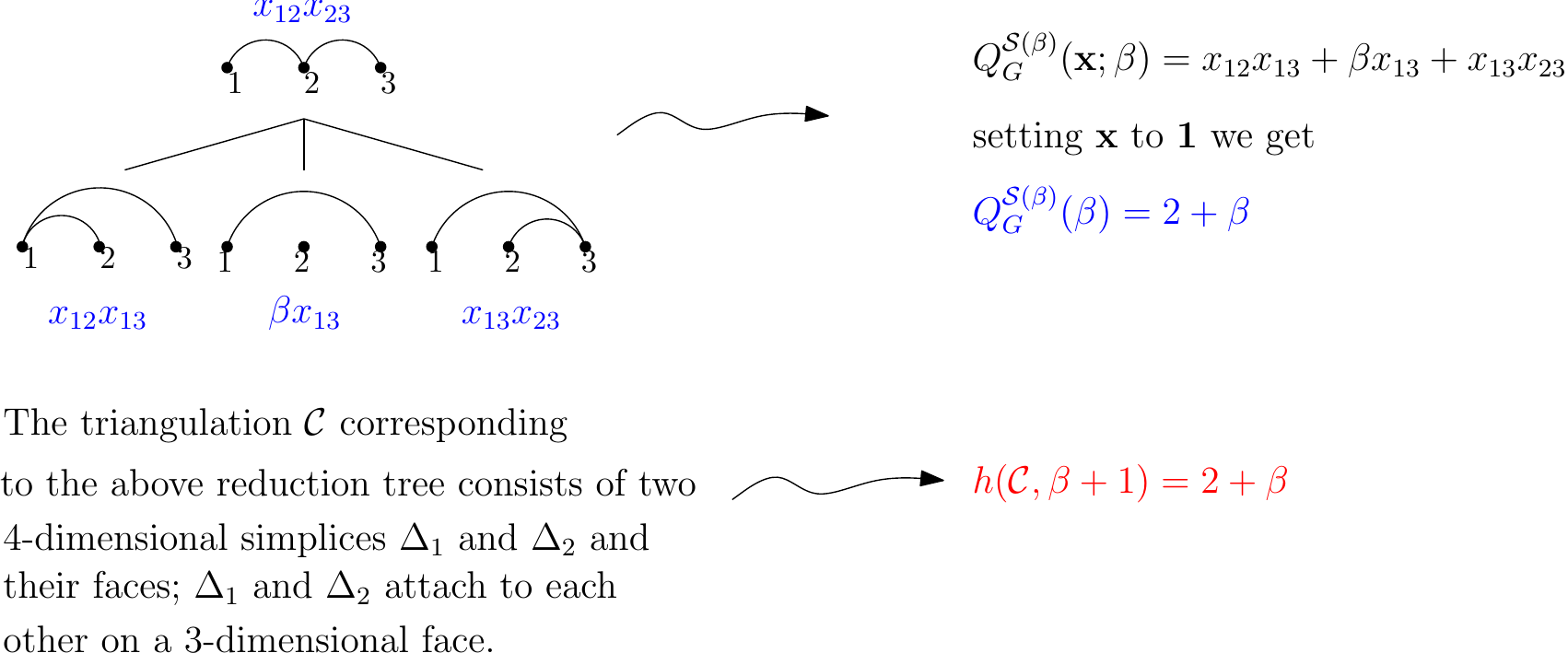}
 \caption{On the left is a reduction tree of $G=([2], \{(1,2),(2,3)\})$ with the monomials corresponding to the graphs noted and a description of the corresponding triangulation $\mathcal{C}$. On the right is the reduced form corresponding to the reduction tree and the (shifted) $h$-polynomial of $\mathcal{C}$. As in Theorem \ref{thm-h} we get $Q_G^{\t}(\b)=h(\C, \b+1)$.}
 \label{h-ex}
 \end{center}
\end{figure}

We illustrate the theorem on an example in Figure \ref{h-ex}. 

To prove Theorem \ref{thm-h}, we use the following lemma, which is implicit in \cite{Stcom}.  For completeness, we include a proof of it. 

\begin{lemma} \label{pure} Let $\C$ be a $(d-1)$-dimensional    simplicial complex homeomorphic to a ball and $f_i^\circ$ be the number of interior faces of $\C$ of dimension $i$. Then 
\begin{equation} \label{h} h(\C, \beta+1)=\sum_{i=0}^{d-1} f_i^\circ \beta^{d-1-i}
\end{equation}
\end{lemma}
\proof We have that 

$$h(\C, \beta+1)=\sum_{F\in \C} (\beta+1)^{\#F}(-\beta)^{d-\#F}=\sum_{i=-1}^{d-1}f_i(\beta+1)^i (-\beta)^{d-i-1}.$$

Thus, \eqref{h} is equivalent to \begin{equation} \label{eq}\sum_{i=-1}^{d-1}f_i(\beta+1)^i (-\beta)^{d-i-1}=\sum_{i=0}^{d-1} f_i^\circ \beta^{d-1-i},\end{equation} which is what we proceed to prove now.
For brevity, denote the left hand side of \eqref{eq} by $LHS$ and the right hand side by $RHS$. Then we see that 

$$[\beta^d]LHS=(-1)^d+\sum_{i=0}^{d-1}f_i(-1)^{d-i-1} {{i+1} \choose 0}.$$ 

For $j \in [d]$ we have that 

$$[\beta^{d-j}]LHS=\sum_{i=j-1}^{d-1}f_i(-1)^{d-i-1} {{i+1} \choose j}.$$ 

On the other hand, $[\beta^d]RHS=0$ and  for $j \in [d]$ we have $[\beta^{d-j}]RHS=f_{j-1}^\circ$. Thus,  \eqref{eq}  is equivalent to

 \begin{equation} \label{eq1} (-1)^d+\sum_{i=0}^{d-1}f_i(-1)^{d-i-1} {{i+1} \choose 0}=0  \end{equation} 
and for $j \in [d]$

 \begin{equation} \label{eq2} \sum_{i=j-1}^{d-1}f_i(-1)^{d-i-1} {{i+1} \choose j}=f_{j-1}^\circ.\end{equation} 

Note that \eqref{eq1} is the familiar Euler characteristic expression of a simplicial complex that is homeomorphic to a  ball. To prove 
\eqref{eq2}  we interpret its  
left hand side   as an inclusion exclusion formula for $f_{j-1}^\circ$. Indeed, $f_{d-1} {{d}\choose {j}}$ is the number of $(j-1)$-faces of the $f_{d-1}$ top dimensional simplices in $\C$, with overcounting. To correct for the overcounting we substract  $f_{d-2} {{d-1}\choose {j}}$, which is the number of $(j-1)$-faces of the $f_{d-2}$ $(d-2)$-simplices of $\C$. Note that here we therefore also take out all $(j-1)$-faces that are on the boundary of $\C$. Again, we have oversubstracted, so we add back $f_{d-3} {{d-2}\choose {j}}$, etc., to arrive to 
$f_{j-1}^\circ$. We used that $\C$ is a pure simplicial complex when we assumed that all $(j-1)$-faces are a face of a $(d-1)$-face. \qed

\bigskip

\noindent \textit{Proof of Theorem \ref{thm-h}.}
Note that, by definition, the reduced form $Q_G^{\t}(\b)$ read off from the reduction tree $R_G^{\O}$, which yields the unimodular triangulation $\C^{\O}$, can be written as

\be \label{uni} Q_G^{\t}(\b)=\sum_{i=1}^{d-1} f^\circ_i \b^{d-1-i},\ee where $\C^{\O}$ is $d-1$ dimensional and $f^\circ_i$ is the number of leaves of $R_G^{\O}$ yielding an $i$-dimensional simplex. The notation $f^\circ_i$ signifies that the union of the open simplices corresponding to the leaves of a reduction tree yield the open polytope $\F_{\widetilde{G}}$.  By Lemma \ref{pure} we have that 
\be \sum_{i=1}^{d-1} f^\circ_i \b^{d-1-i}=h(\C^{\O}, \b+1).\ee Since by Lemma \ref{leaves} we know that $Q_G^{\t}(\b)$ does not depend on the particular reduction tree, we are done. \qed

\section{A description of the leaves of $R_G^{\O}$}
\label{sec:leaves}

In this section we  describe all the leaves of the reduction tree $R_G^{\O}$ in terms of shellings of the triangulation obtained from $R_G^{\O}$ as in Theorem \ref{tri}. 
Since we know that these triangulations are regular, it follows that they are also shellable.

\bt \label{thm-shell2}  Let $\F_{\widetilde{F_1}}, \ldots, \F_{\widetilde{F_l}}$ be a shelling order of the simplicial complex arising from $R_G^{\O}$. Let $$P_i:=\{\{Q_1^i, \ldots, Q_{f(i)}^i\}\}=\{\{F_i \cap F_j \mid 1\leq j<i, |E(F_i \cap F_j)|=|E(F_i)|-1\}\}.$$  Then 
\be \label{formal} \sum_{i=1}^l \prod_{j=1}^{f(i)} (F_i+Q_j^i)\ee

\noindent is the formal sum of the  set of the leaves of $R_G^{\O}$, where the product of graphs is their intersection. If $f(i)=0$ we define  $\prod_{j=1}^{f(i)} (F_i+Q_j^i)=F_i$.

\et

Before proving Theorem \ref{thm-shell2}, we record  a few properties of flow polytopes that easily follow from the above considerations and from the fact that the dimension of $\F_G$ is $|E(G)|-|V(G)|+1$ \cite{BV2}. In both lemmas the meanings of 
$\F_{\widetilde{H}}$ for a node $H$ of $R_G^{\O}$ or intersection of such nodes is as in Definition \ref{leaf}, which is the key to the proofs that are left to the interested reader.

\bl \label{coord} Let $G_1$ and $G_2$ be two leaves of $R_G^{\O}$. Then 

\be \label{int}  \F_{\widetilde{G_1}}\cap \F_{\widetilde{G_2}}=\F_{\widetilde{G_1\cap G_2}}.\ee Moreover, $G_1\cap G_2$ is a leaf of $R_G^{\O}$ if and only if $\F_{\widetilde{G_1}}\cap \F_{\widetilde{G_2}}$ is not contained in the boundary of $\F_{\widetilde{G}}$. 
\el

\bl \label{dimension}  Let $G_1$ and $G_2$ be two leaves of $R_G^{\O}$. The dimension of  $\F_{\widetilde{G_1\cap G_2}}$ is $|E(G_1\cap G_2)|+|V(G_1\cap G_2)|-1$.
\el



\noindent \textit{Proof of Theorem \ref{thm-shell2}.}  By Lemmas \ref{coord} and \ref{dimension} we see that if $\F_{\widetilde{F_1}}, \ldots, \F_{\widetilde{F_l}}$ is a shelling order, then the set of facets on which $\F_{\widetilde{F_i}}$ attaches to $\F_{\widetilde{F_1}}, \ldots, \F_{\widetilde{F_{i-1}}}$ is $\{\F_{\widetilde{Q}} \mid Q \in P_i\}$. Moreover, since the intersection of two top dimensional simplices of a triangulation of  $\F_{\widetilde{G}}$ is not contained in the boundary of  $\F_{\widetilde{G}}$, it follows that every element of $P_i$, $i \in [l]$, appears in $R_G^{\O}$ by the second part of Lemma \ref{coord} (and it is a leaf since it is the intersection of two alternating graphs, so it is alternating itself).  Using this same argument repeatedly and the fact that we can built up the polytope piece by piece by following the shelling, we obtain Theorem \ref{thm-shell2}. 
\old{
By Theorem \ref{tri} the vertices of the simplices corresponding to the leaves of $R_G^{\O}$ are mutually coherent in the $\widetilde{\O}$ framing of $\widetilde{G}$. Using this  it is easy to see that for any two leaves $G_1$ and $G_2$ of  $R_G^{\O}$ we have that  $$\F_{\widetilde{G_1}}\cap \F_{\widetilde{G_2}}=\F_{\widetilde{G_1\cap G_2}}.$$ Moreover, $G_1\cap G_2$ is a leaf of $R_G^{\O}$ exactly if $\F_{\widetilde{G_1}}\cap \F_{\widetilde{G_2}}$ is not contained in the boundary of $\F_{\widetilde{G}}$. Using this and how the polytope is built up piece by piece in the shelling Theorem \ref{thm-shell2} follows.}
\qed
\old{

In order to apply Theorem \ref{thm-shell2}, it would be nice to have an explicit shelling order on hand. The next theorem, which is proved in \cite{h-poly2} accomplishes this. See Figure \ref{fig:o} for an illustration of Theorem \ref{thm-shell2} with the shelling order from Theorem \ref{thm-shell}.

Let  $F_1, \ldots, F_l$  be  the full-dimensional leaves of $R_G^{\O}$ ordered by depth-first search order. Remember that we have an embedding of $R_G^{\O}$ in the plane where $G$ is the root and the graphs $G_1, G_2, G_3$ as in \eqref{graphs1} are the left, right, middle child, respectively.

\bt \label{thm-shell}  \cite[Theorem ]{h-poly2}$\F_{\widetilde{F_1}}, \ldots, \F_{\widetilde{F_l}}$ is a shelling order of the triangulation of $\F_{\widetilde{G}}$ from Theorem \ref{tri}. \et 

 }

\section{Nonnegativity results about  reduced forms}
\label{sec:nonneg}

This section is devoted to two nonnegativity results, which are consequences of the above considerations.  

\bt \label{pos} The polynomial $Q_G^{\t}(\beta-1)$ is  a polynomial in $\b$ with nonnegative coefficients.
\et

\proof Recall that by Theorem \ref{thm-h} we have that 
\be Q_G^{\t}(\b-1)=h(\C, \b),  \ee where $\C$ is a unimodular triangulation of    $\F_{\widetilde{G}}$. Let $\C$ be the abstract  simplicial complex obtained from $R_G^{\O}$, as in Theorem \ref{tri}. Since by Theorem \ref{tri} this triangulation is regular, and therefore it is shellable,  we get  that  $h_i$ is equal to the number of top dimensional simplices which attach on $i$ facets to the union of previous simplices in a shelling order. 
\qed

\medskip

Using Theorem \ref{pos} we are ready to prove a special case of Kirillov's  \cite[Conjecture 2]{k2}.

\begin{conjecture} \cite[Conjecture 2]{k2} Let $k_1, \ldots, k_{n-1}$ be a sequence of nonnegative integers and let $M=x_{12}^{k_1}x_{23}^{k_2}\cdots x_{n-1, n}^{k_{n-1}}$. Then the reduced form of $M$    evaluated at $\x=(1, \ldots, 1)$ and $\b-1$  in  $\widetilde{ACYB_n}(\beta)$ is a polynomial in $\beta$ with nonnegative coefficients. 
\end{conjecture}

\bt \label{special} \cite[cf. Conjecture 2]{k2}. The reduced form of $x_{12}x_{23}\cdots x_{n-1,n}$ evaluated at $\x=(1, \ldots, 1)$ and $\b-1$ in  $\widetilde{ACYB_n}(\beta)$  is a polynomial in $\b$ with nonnegative coefficients.
\et

\proof  In \cite{m2} it is proved that the monomial $x_{12}x_{23}\cdots x_{n-1,n}$  can be reduced in $\widetilde{ACYB_n}(\beta)$ so that the resulting monomials in the reduced form have no variables of the form $x_{ij}x_{jk}$, $i<j<k$. Thus, reduced form of $x_{12}x_{23}\cdots x_{n-1,n}$ evaluated at $\x=(1, \ldots, 1)$ and $\b-1$ in  $\widetilde{ACYB_n}(\beta)$  is equal to the reduced form of $x_{12}x_{23}\cdots  x_{n-1,n}$ evaluated at $\x=(1, \ldots, 1)$ and $\b-1$ in  $\t$.  Therefore we can apply Theorem \ref{pos} to obtain this result.
\qed

\section{Reduced forms and Ehrhart series}
\label{sec:ehr}

In this section we connect the $h$-polynomials of triangulations of flow polytopes  to the  Ehrhart series of flow polytopes, and using Theorem \ref{thm-h} we tie this in with reduced forms in the subdivision algebra. As a corollary to our results, we generalize \cite[Theorem 3.10]{k2}.

Recall that for a  polytope $\mathcal{P}\subset \mathbb{R}^{N}$, the {$t^{th}$ dilate} of $\mathcal{P}$ is $\displaystyle t \mathcal{P}=\{(tx_1, \ldots, tx_{N}) \mid  (x_1, \ldots, x_{N}) \in \mathcal{P}\}.$ The number of lattice points of $t\mathcal{P}$, where $t$ is a nonnegative integer and $\mathcal{P}$ is a convex polytope, is given by the {\bf Ehrhart function} $i({\mathcal{P}}, t)$. If $\mathcal{P}$ has integral vertices then $i({\mathcal{P}}, t)$ is a  polynomial.

The following is a well-known result specialized to  $\F_{\widetilde{G}}$.

\begin{lemma} \label{thm-ehr} Let  $\C$ be a unimodular triangulation of  $\F_{\widetilde{G}}$.
Then
\be \label{h-ehrhart1} h(\C, \b)=\sum_{m\geq 0}(i(\F_{\widetilde{G}}, m)\b^m)(1-\b)^{\dim(\F_{\widetilde{G}})+1}.\ee
\end{lemma}

\bc \label{cor} We have 
\be Q_G^{\t}(\b-1)=\sum_{m\geq 0}(\l(\F_{\widetilde{G}}, m)\b^m)(1-\b)^{\dim(\F_{\widetilde{G}})+1}.\ee  \ec

\proof Follows directly from Theorem \ref{thm-h} and Lemma  \ref{thm-ehr}.\qed

\bc{\cite[Theorem 3.10]{k2}} 
\be Q_{K_n}^{\t}(\b-1)=\sum_{m\geq 0}(\l(CRY_{n+1}, m)\b^m)(1-\b)^{{n+1 \choose 2}}.\ee  
\ec

\proof Follows from Corollary \ref{cor}  for $G=K_n$, since $\dim (\F_{\widetilde{K_n}})={{n+1 \choose 2}-1}$ and  $\l(CRY_{n+1}, m)=i(\F_{\widetilde{K_n}}, m)$, as explained in \cite{m-prod}. \qed

To state our final result, also a corollary of Corollary \ref{cor}, relating the reduced forms to Kostant partition functions, we remind the reader of the following definition.

 The {\bf Kostant partition function}  $K_G$ evaluated at the vector $\v \in \Z^{n+1}$ is defined as

\begin{equation} \label{kost} K_G(\v)= \#  \{ (b_{k})_{k \in [N]} \mid \sum_{k \in [N]} b_{k}  \a_k =\v \textrm{ and } b_{k} \in \Z_{\geq 0} \},\end{equation}

\noindent where $[N]=\{1,2,\ldots,N\}$ and  $\{\{\a_1, \ldots, \a_N\}\}$ is the multiset of vectors corresponding to the multiset of edges of $G$ under the correspondence which associates an edge  $(i, j)$, $i<j$, of $G$ with a positive type $A_n$ root $e_i-e_j$, where $e_i$ is the $i$th standard basis vector in $\mathbb{R}^{n+1}$.
In other words, $K_G(\v)$ is the number of ways to write the vector $\v$ as a $\mathbb{N}$-linear combination of the positive type $A_n$ roots (with possible multiplicities) corresponding to the edges of $G$, without regard to order.

\begin{corollary} For any graph $G$ we have 
\be Q_G^{\t}(\b-1)=\sum_{m\geq 0}(K_{\widetilde{G}}(m, 0, \ldots, 0, -m)\b^m)(1-\b)^{\#E(G)+\#V(G)},\ee where $K_{\widetilde{G}}$ is the Kostant partition function for ${\widetilde{G}}$. \end{corollary}

\proof  Follows from Corollary \ref{cor}   since $\l(\F_{\widetilde{G}}, m)=K_{\widetilde{G}}(m, 0, \ldots, 0, -m)$ is a simple corollary of the definitions of these objects and $\dim(\F_{\widetilde{G}})=\#E(\widetilde{G})-\#V(\widetilde{G})+1=\#E(G)+\#V(G)-1.$ For  detailed explanations of both of these and related results see \cite{mm}. \qed


  \section*{Acknowledgement}
  The author  thanks Federico Ardila, Louis Billera and Ed Swartz for several  interesting  discussions regarding this work. The author also thanks Alejandro H. Morales for many inspiring conversations about flow polytopes over the course of the years.

\bibliography{biblio-kir}
\bibliographystyle{plain}

 \end{document}